\documentclass[pdflatex,sn-mathphys-num]{sn-jnl}

\usepackage{graphicx}%
\usepackage{multirow}%
\usepackage{amsmath,amssymb,amsfonts}%
\usepackage{amsthm}%
\usepackage{mathrsfs}%
\usepackage[title]{appendix}%
\usepackage{xcolor}%
\usepackage{textcomp}%
\usepackage{manyfoot}%
\usepackage{booktabs}%
\usepackage{algorithm}%
\usepackage{algorithmicx}%
\usepackage{algpseudocode}%
\usepackage{listings}%

\theoremstyle{thmstyleone}%
\theoremstyle{thmstyletwo}%
\newtheorem{example}{Example}%
\newtheorem{remark}{Remark}%

\theoremstyle{thmstylethree}%

\usepackage{subfigure}
\usepackage{caption}

\usepackage{adjustbox} %

\begin{document}

\title[Article Title]{A persistent-homology-based Bayesian prior for potential coefficient reconstruction in an elliptic PDE}


\author*[1]{\fnm{Zhiliang} \sur{Deng}}\email{dengzhl@uestc.edu.cn}

\author[1]{\fnm{Haiyang} \sur{Liu}}\email{202321110210@std.uestc.edu.cn}

\author[2]{\fnm{Xiaofei} \sur{Guan}}\email{guanxf@tongji.edu.cn}
\author[2]{\fnm{Zhiyuan} \sur{Wang}}\email{2211187@tongji.edu.cn}
\author*[3]{\fnm{Xiaomei} \sur{Yang}}\email{yangxiaomath@swjtu.edu.cn}

\affil*[1]{\orgdiv{School of Mathematical Sciences}, \orgname{University of Electronic Science and Technology of China}, \orgaddress{\street{Xiyuan Ave}, \city{Chengdu}, \postcode{611731}, \state{Sichuan}, \country{China}}}

\affil[2]{\orgdiv{School of Mathematical Sciences}, \orgname{Tongji University}, \orgaddress{\street{Siping Road}, \city{Shanghai}, \postcode{200092}, \state{Shanghai}, \country{China}}}

\affil[3]{\orgdiv{School of Mathematics}, \orgname{Southwest Jiaotong University}, \orgaddress{\street{Xipu Campus}, \city{Chengdu}, \postcode{610097}, \state{Sichuan}, \country{China}}}


\abstract{

We address the reconstruction of a potential coefficient in an elliptic partial differential equation from distributed observations within the Bayesian framework. The choice of prior distribution is crucial in such inverse problems, particularly when the target function exhibits sharp discontinuities that conventional Gaussian priors fail to capture effectively. To overcome this limitation, we introduce a novel prior based on persistent homology (PH), which quantifies and encodes the topological features of candidate functions through their persistent pairs. To ensure a well-defined distribution in infinite-dimensional spaces, the prior is constructed with respect to a Gaussian reference measure. A significant advantage over classical approaches is that the PH prior only requires the unknown functions to belong to a suitable topological space, which substantially enhances its applicability. Numerical results demonstrate that the proposed PH prior outperforms the Gaussian prior and achieves a modest yet consistent improvement over the classical total variation (TV) prior.

 }

\keywords{ Inverse potential problems, topological prior, Bayesian approach, persistent homology, TV prior}



\maketitle

\section{Introduction}\label{sec1}
The inverse potential problem  in the paper is governed by the elliptic equation
\begin{align}\label{eqn1.1}
\left\{
\begin{aligned}
&-\Delta u+qu=f &\text{in} && \Omega, \\
& u=0 &\text{on}  &&\partial \Omega,
\end{aligned}
\right.
\end{align}
where $\Omega\subset\mathbb{R}^{\mathfrak{s}}$ ($\mathfrak{s}=1, 2$) is a bounded open domain, 
and the function $f$ is the known source term. The potential $q$ belongs to the admissible set $K$ defined as
\begin{align}\label{eqnK}
K=\{q\in L^\infty(\Omega): c_0\leq q(x)\leq c_1 \,\, \text{a.e.}\,\, \text{in}\,\, \Omega \},
\end{align} 
with $0\leq c_0<c_1<\infty$. We collect the observational data of the solution $u(x)$ to \eqref{eqn1.1} on $\Omega$ by 
\begin{align}\label{eqn1.2}
u^\eta(x)=u(x)+\eta(x),\,\, x\in \Omega,
\end{align}
where $\eta$ denotes the measurement noise. Our aim is to estimate the potential function $q$ from the noisy observation $u^\eta$. 

The relative issues occur in various applied fields. 
Choulli's work \cite{Choulli2021}  notes that the inverse potential problem is  connected to dynamic elastography, microwave imaging, and acousto-optic imaging.
Bal et al \cite{Bal2010} apply it to  the photoacoustic tomography problem.
In some real-world cases, the inverse potential problem is also inherently linked to time-dependent phenomena, e.g., the reconstruction of the radiativity coefficient in heat equation \cite{Bal2010, Jin2023, Pennes1948, Trucu2010, Yamamoto2001} and the 
reconstruction of  the potential in Schr\"{o}dinger equation \cite{Baudouin2002, Baudouin2008}.   Due to the broad applicability, inverse potential problems have attracted considerable research attention.
The related theoretical results have been established by some authors.  Bal et al \cite{Bal2010} show that the uniqueness and the stability of the reconstruction for an inverse diffusion problem with internal data.
Jin et al \cite{Jin2023} establish a weighted conditional stability estimates under  mild conditions on the problem data. In \cite{Choulli2021}, the author derives both Lipschitz and H\"{o}lder stability estimates. These works provide sufficient conditions under which the solution is uniquely determined and remains stable with respect to perturbations in the data.
For convergence analysis and uncertainty quantification, Nickl \cite{Nickl2020_a} discusses posterior contraction rates and proves a Bernstein–von Mises theorem, ensuring consistency of the Bayesian posterior and providing rigorous uncertainty quantification.  In \cite{Nickl2020_b}, the authors give the convergence rate of the maximum a posteriori (MAP) estimator within the Bayesian framework.
Jin et al. \cite{Jin2023} derive a convergence rate in a weighted $L^2(\Omega)$ using Tikhonov regularization with an $H^1$-seminorm penalty. And they  \cite{Jin2025} further investigates the convergence of the regularized solution and its finite element approximation, providing a stochastic $L^2(\Omega)$ convergence rate with high probability.


In these literature,  the Bayesian approach is a main technique in inverse potential problems. 
It offers a flexible and probabilistically grounded framework for incorporating uncertainty in both data and model parameters. Moreover, the Bayesian posterior distribution serves as a rigorous basis for uncertainty quantification, probabilistic inference, and informed decision making.
 The theoretical developments of Bayesian inversion have fueled the practical success  in a wide range of inverse problems  \cite{Beskos2017, Bui-Thanh2014,  Carpio2020, Dashti2017,  Furuya2024, Giordano2020, Huang2021, Kaipio2005,  Kaipio2019, Laghrib2022,  Li2020,  Stuart2010, Tarantola2005}. 
One key in this approach is the specification of the prior distribution, which encodes prior beliefs or structural assumptions about the unknowns.
Considerable efforts have been dedicated to designing appropriate priors. The common  priors are derived from Markov random fields  (e.g. \cite{Bardsley2013, Iglesias2016}),  which offer a rich modeling framework.
A variety of other priors have also been developed to address specific problems.
 In \cite{Nickl2020_a}, Nickl introduces a wavelet-based prior for the reconstruction of potential coefficients.
To address inverse problems with less regular parameters, Wang et al \cite{Wang2017} present an $l_1$-type prior that includes both the TV prior and the Besov space $B_{1,1}^s$ prior. In \cite{Yao2016}, a hybrid TV-Gaussian prior is developed to better capture sharp discontinuities in the target function.
The incorporation of geometric information into the prior is another fruitful line of research, especially in problems involving geometric reconstructions. Iglesias et al \cite{Iglesias2014} give a geometry-informed prior formulation, where prior assumptions are guided by the structural characteristics of the unknown domain. In  \cite{Carpio2020, Carpio2023},  Carpio et al propose a novel prior construction strategy grounded in the topological energy landscape of objective functions, which they successfully apply to inverse scattering and wave equation parameter estimation.

In this paper, we introduce a new prior based on the topological tool of persistent homology (PH-based prior) and apply it to potential reconstruction within a Bayesian framework. This prior differs from those proposed by Carpio et al. \cite{Carpio2020}. It aims to capture the topological variation of the unknown and explore its inherent structural information. Through these topological features, latent shape information can be extracted and analyzed. This provides a fresh perspective for understanding the unknown functions.
Our approach with the prior achieves substantially better numerical performance than Gaussian priors in estimating non-smooth functions, particularly in cases involving abrupt changes (e.g., sharp jumps or step discontinuities). 
In formulation, the PH-based prior is similar to the TV-based prior \cite{Yao2016}.  In contrast, our prior can be defined on a more general topological space and preserves richer topological information than the TV-based prior. This makes it applicable to a wider range of fields  and particularly useful for complex inverse problems. In this sense,  the prior is an extension of the TV-based prior in \cite{Yao2016}.

This paper is organized as follows: In Section 2, we present the necessary background on persistent homology;  Section 3 introduces the Bayesian approach to the inverse potential problem. 
Building upon this foundation, we propose a PH-based prior in Section 4 and apply it to the inverse potential problem. To enable efficient sampling, we decompose the PH-based prior into two components and develop a delayed acceptance algorithm in Section 5. Finally, we present our conclusions in Section 6.

\section{The persistent homology}

This section introduces some knowledge on persistence diagram and persistent homology which form a concrete basement for our prior. For more details one can refer to \cite{Edelsbrunner2014, Kaczynski2004, Zheng2015, Zomorodian2005}.

For a function $q:\Omega\to\mathbb{R}$, 
we define the sublevel set for each $c\in\mathbb{R}$ by
$$
\mathcal{S}_c = \{x\in\Omega \mid q(x)\le c\}.
$$
The topology of $\mathcal{S}_c$ changes only when $c$ passes through a critical value of $q$. Its homology groups $H_k(\mathcal{S}_c)$, $k=0,1,\dots$, capture these topological features, where $k$ indicates the dimension (connected components for $k=0$, loops for $k=1$, voids for $k=2$, etc.).

In practice, we compute these features by discretizing $\Omega$ into a simplicial complex $\mathcal{K}$, based on a partition $\{x_i\}$. Each partition node $x_i$ is identified with a vertex $v_i\in\mathcal{K}$, so that the vertex set $\{v_i\}$ corresponds exactly to the partition nodes. The function $q$ is then defined at vertices by $q(v_i) := q(x_i)$; for simplicity, we identify this vertex function with $q$ itself. For a simplex $\sigma = [v_0,\dots,v_m]$, we extend $q$ by
$$
q(\sigma) = \max\{q(v_0),\dots,q(v_m)\},
$$
which  naturally induces the subcomplex
$$
\mathcal{K}_c = \{\sigma\in \mathcal{K}\mid q(\sigma)\le c\}.
$$
Clearly,  $\mathcal{K}_{c^-}\subseteq \mathcal{K}_c$ for $c^{-}\leq c$;  this nested structure gives rise to the filtration  $\{\mathcal{K}_c\}$.
By the Nerve theorem and simplicial approximation \cite{Edelsbrunner2014}, the homology groups $H_k(\mathcal{K}_c)$ and $H_k(\mathcal{S}_c)$ are isomorphic:
$$
H_k(\mathcal{K}_c)\cong H_k(\mathcal{S}_c).
$$
Persistent homology tracks the birth and death of homology classes as $c$ increases. Algorithmically, this is realized by pairing simplices in the filtration $\{\mathcal{K}_c\}$. When a simplex $\sigma$ enters $\mathcal{K}_c$, it may create a new homology class or destroy an existing one. 
Denoting the homology group at the previous step by $H_k(\mathcal{K}_{c^-})$, the change in the homology group can be expressed as
\begin{align}
\begin{aligned}
&H_k(\mathcal{K}_c)=H_k(\mathcal{K}_{c^-})\oplus \langle\bar{\sigma}\rangle,\quad \text{birth of a new }k\text{-class}\,\, \bar{\sigma};\\
&H_k(\mathcal{K}_c)=H_k(\mathcal{K}_{c^-})/ \langle\bar{\tau}\rangle,\quad \text{death of an existing }k\text{-class}\,\, \bar{\tau}.
\end{aligned}
\end{align}
If the simplex neither creates nor destroys a class, then $H_k(\mathcal{K}_c) \cong H_k(\mathcal{K}_{c^-})$.
This  provides an algorithmic way to track the evolution of homology groups and enables the computation of persistent homology. 
The essence of the topology change in the filtration $\{\mathcal{K}_c\}$  is precisely reflected by persistence pairs: if a class is born at $\sigma$ and later destroyed at $\tau$, we record the persistence pair $(\sigma, \tau)$; if the class never dies, it is paired with $\infty$.
The persistence of a class is computed as 
$|q(\tau)-q(\sigma)|$, the difference between its death time $q(\tau)$
 and birth time $q(\sigma)$.  Therefore,  the topology of the function   $q:\Omega\rightarrow \mathbb{R}$ can be effectively characterized by computing persistence pairs $(\sigma, \tau)$ and tracking lifetimes $|q(\tau)-q(\sigma)|$ of homology classes.

 


In one dimension, $\Omega\subset\mathbb{R}$ is discretized into ordered grid points $\{x_i\}$, and $\mathcal{K}$ consists of vertices $v_i$ (corresponding to $x_i$) and edges $[v_i,v_{i+1}]$. The function $q$ is represented at vertices by $q(v_i)=q(x_i)$, and extended to edges via
\begin{align}\label{maxv}
q((v_i,v_{i+1})) = \max\{q(v_i),q(v_{i+1})\}.
\end{align}
Here, its topological features only include connected components. Each local minimum creates a new component, which is destroyed when it merges with an older one. Thus, 1D persistent homology (edge induced filtration) reduces to a vertex pairing scheme: a vertex $v_k$ corresponding to a local minimum creates a component, a vertex $v_l$ marks its destruction, and the persistence is
$$
\mathrm{pers}(v_k,v_l) = |q(v_l)-q(v_k)|.
$$
For components that never merge, the corresponding vertex is paired with $\infty$, and their persistence is taken as
$$
\mathrm{pers}(v_k,\infty) = \infty.
$$
We collect all  pairs with finite persistence in $P_+$, and similarly consider $-q$ to obtain $P_-$. The persistence distance \cite{Zheng2015} is then 
$$
\|q\|_{\rm per} := \sum_{(v_k,v_l)\in P_+}\mathrm{pers}(v_k,v_l) + \sum_{(v_k,v_l)\in P_-}\mathrm{pers}(v_k,v_l).
$$
It can be seen that $\|q\|_{\rm per}$ encodes the total “ups and downs” of $q$ excluding the overall global range. 
Let $x_{i_0}<\dots<x_{i_r}$ denote all local extrema of $q$. Then the discrete total variation satisfies
$$
\mathrm{TV}(q) = \sum_{j=0}^{J-1}|q(x_{j+1})-q(x_j)| = \sum_{t=0}^{r-1}|q(x_{i_{t+1}})-q(x_{i_t})|.
$$
The  persistence pairs in $P_+$ and $P_-$ account for all adjacent-extrema differences except the overall range $\max_j q(v_j)-\min_j q(v_j)$. Hence, we have the identity
$$
\|q\|_{\rm per} + \max_{i, j} |q(v_i)-q(v_j)| = \mathrm{TV}(q),
$$
which establishes a simple and transparent link between  $1$D persistence and classical discrete total variation.

\begin{remark}
In \cite{Zheng2015}, the extension of $q$ from vertices to simplices is realized via the piecewise linear spline interpolant.  We use \eqref{maxv} to extend $q$ from vertices to simplices. 
It induces the same discrete filtration for $0$D homology as in \cite{Zheng2015}. Therefore,  both approaches lead to the same vertex pairing and persistence values.  However the maximum-extension rule is simpler and directly applicable to functions $q \in L^\infty(\Omega)$ that may be discontinuous.
\end{remark}

\section{The Bayesian approach to the potential coefficient reconstruction}

For simplicity, we assume that the measurements are collected at $m$ uniformly distributed points $\{x_i\}_{i=1}^m \subset \Omega$. The inverse problem is then formulated as
$$
u^\eta(x_i) = G(q)(x_i) + \eta(x_i), \quad i=1,\dots,m,
$$
where $G:K \to \mathbb{R}^m$ is the forward operator and $\eta \sim \mathcal{N}(0,\Sigma)$ denotes Gaussian noise.
Under this assumption, the likelihood of $u^\eta$ given  $q$ is
\begin{align}
\pi(u^\eta|q)\propto \exp\left(-\Phi(q; u^\eta)\right),
\end{align}
where the data misfit functional
\begin{align}
\Phi(q; u^\eta):=\frac{1}{2}\left\|G(q)-u^\eta\right\|^2_{\Sigma}=\frac{1}{2}\left|\Sigma^{-1/2}\left(G(q)-u^\eta\right)\right|^2.
\end{align}
This term corresponds to  the data fidelity term in deterministic inverse problems. The Bayesian approach provides a principled framework for inference under uncertainty, where the prior distribution encodes prior knowledge about the unknowns and quantifies uncertainty before incorporating the data.
 Here we assume that the prior measure of \(q\) is  denoted by \(\mu_{\rm pr}\). While $\mu_{\rm pr}$ is  defined on $L^\infty(\Omega)$, the dual space of $L^\infty(\Omega)$ is notoriously complicated.  Therefore, we restrict our consideration to the corresponding measures on $L^2(\Omega)$.
 When the posterior measure \(\mu_{\text{post}}\) is absolutely continuous with respect to the prior measure \(\mu_{\rm pr}\), it admits the Radon-Nikodym derivative
\begin{align}
\frac{\rm{d}\mu_{\rm post}}{\rm{d}\mu_{\rm pr}}(q)=\frac{1}{Z}\exp\left(-\Phi(q; u^\eta)\right),
\end{align}
where $Z$ is a normalization constant. 

The Gaussian prior, $\mu_{\rm pr} = \mu_0$, where $\mu_0 = N(0, C_0)$, is the most commonly used prior in Bayesian inverse problems. In practical terms, such a prior reflects the assumption that the unknown parameter $q$ has a zero mean, with the covariance $C_0$ encoding our beliefs about the uncertainty or spread of $q$. 

For convenience, we begin by introducing the Gaussian measure on $L^2(\Omega)$, which serves as our reference measure.
For every \( \iota \in L^2(\Omega)^* \), if \( \mu_0 \circ \iota^{-1} \) is a Gaussian measure on \( \mathbb{R} \), then \( \mu_0 \) is defined as a Gaussian measure on \( L^2(\Omega) \).  
The covariance operator \( C_0 \) associated with \( \mu_0 \) is given by:  
\begin{align}\label{eqn16}
C_0(\iota_1, \iota_2) = \int_{L^2(\Omega)} \iota_1(q)\iota_2(q) \, d\mu_0(q) - \mathbb{E}[\iota_1] \mathbb{E}[\iota_2], \quad \iota_1, \iota_2 \in L^2(\Omega)^*,
\end{align}  
where \( \mathbb{E}[\iota] \) denotes the mean of \( \mu_0 \), defined as:  
\begin{align}
\mathbb{E}[\iota] = \int_{L^2(\Omega)} \iota(q) \, d\mu_0(q), \quad \text{for every } \iota \in L^2(\Omega)^*.
\end{align}  
Next we interpret \eqref{eqn16} in a more intuitive manner and express the covariance operator differently. According to the Hahn-Banach theorem, we know that \( \iota(q) \) serves, to some extent, as a means of extracting the coordinate of \( q \). This implies that \( C_0(\iota_1, \iota_2) \) represents the covariance between different coordinates, $\iota_1, \iota_2$,  of \( q \). By the Riesz representation theorem, we no longer distinguish between \( L^2(\Omega) \) and its dual \( L^2(\Omega)^* \) in the following. Let \( \{\varphi_n\}_{n=1}^\infty \) be an orthonormal basis of \( L^2(\Omega) \). The Fourier expansion of $q$ is given by $q(x)=\sum_{n=1}^\infty q_n \varphi_n(x)$. 
As taking $\iota$ as $\varphi_i$, $\varphi_j$, it gets the coordinates $q_i$, $q_j$ of $q$ by $\iota(q)$. It can be seen that the mean and covariance operators are bounded linear
and bilinear functional respectively. The bounded bilinear functional determines a bounded linear operator as follows:  
\begin{align}  
C_0(\iota_1, \iota_2) = \langle C_0 \iota_1, \iota_2 \rangle,  
\end{align}  
where the same notation is used for both the bilinear functional and the operator. In the present paper, we consider two commonly used forms of \( C_0 \). 
The first is to define the covariance operator \( C_0 \) as an integral operator with a mean squared exponential kernel:  
\begin{align} \label{prior1} 
C_0^{\exp} \iota = \int_{\Omega} \exp\left(-\frac{|x-y|^2}{2l^2}\right) \iota(y) \, dy,
\end{align}  
where $l>0$ is a length-scale parameter.
In some cases of periodic structure, we usually use the periodic version of $C_0^{\exp}$, i.e., the periodic mean squared exponential kernel covariance operator $C_0^{\rm per}$
\begin{align} \label{prior2}
C_0^{\rm per} \iota = \int_{\Omega} \exp\left(-\frac{2\sin^2(\pi|x-y|/p)}{l^2}\right) \iota(y) \, dy, 
\end{align}  
where $l>0$ is a length-scale parameter and $p>0$ is a periodicity parameter. 
The second is to define
 \begin{align} \label{prior3}
 C_0^{-\Delta} = (-\Delta)^{-s}, \,\, s>0
 \end{align}
 with Dirichlet boundary conditions.  

Gaussian priors with covariance operators, as described in equations \eqref{prior1}-\eqref{prior3}, promote smoothness by incorporating prior information that encourages regularity in the underlying function. Functions sampled from these priors are typically smooth. 
However, non-smooth unknown functions $q$ often arise in practical scenarios, such as in the infinite potential well model in the Schrödinger equation.
In such cases, Gaussian priors may struggle to capture the non-smooth characteristics of the function. To address this limitation, we introduce an additional constraint for the unknown function $q$ within the framework of the Gaussian prior. This results in a hybrid prior form, which combines the strengths of both the probabilistic Gaussian framework and regularization techniques. The concept of hybrid priors originates from the work of Yao et al. \cite{Yao2016}, who proposed a more flexible prior that incorporates both uncertainty modeling and total variation methods. This hybrid approach enhances the modeling capabilities in complex systems, where both uncertainty and structural information are critical for accurate representation.

\section{The  PH-Gaussian prior}

We start by reviewing the work of Yao et al.~\cite{Yao2016}.
Rather than simply setting \(\mu_{\rm pr} = \mu_0\), the prior measure is defined as
\[
\frac{{\rm d}\mu_{\rm pr}}{{\rm d}\mu_0}(q) \propto \exp\left(-R(q)\right),
\]
where \(R(q)\) incorporates additional information beyond the Gaussian prior. Under this assumption, it immediately follows that the Radon-Nikodym derivative of \(\mu_{\rm post}\) with respect to \(\mu_0\) is given by
\begin{align}\label{eqn18}
\frac{{\rm d}\mu_{\rm post}}{{\rm d}\mu_0}(q) \propto \exp\left(-\Phi(q; u^\eta) - R(q)\right):=\exp\left(-V(q)\right),
\end{align}
which recovers the standard formulation with Gaussian priors. 
Let the sample space be the Sobolev space \(H^1(\Omega)\), defined as follows:
\[
H^1(\Omega) = \left\{ q \in L^2(\Omega) \mid \partial_x^\alpha q \in L^2(\Omega) \,\, \text{for all} \,\, |\alpha| \leq 1 \right\},
\]
where $\alpha=(\alpha_1, \alpha_2, \cdots, \alpha_{\mathfrak{s}})$ and $|\alpha|=\sum_{i=1}^{\mathfrak{s}}\alpha_i$, and the associated norm is 
\begin{align*}
\|q\|_{H^1}=\sum_{|\alpha|\leq 1}\|\partial_x^\alpha q\|_{L^2(\Omega)}.
\end{align*} 
The regularization term is chosen to be the TV seminorm \cite{Yao2016}
\begin{align}\label{eqntv}
R_{\rm TV}(q)=\lambda\|q\|_{\rm TV}=\lambda\int_{\Omega} |\nabla q| dx,
\end{align}
where $\lambda$ is a prescribed positive constant.  Moreover, the posterior can be written as
\begin{align}\label{eqntv_c}
\frac{d\mu_{\rm post}}{d\mu_0}=\exp\left(-\Phi(q; u^\eta)-R_{\rm TV}(q)\right).
\end{align}

We next introduce a novel hybrid prior that combines the strengths of Gaussian distributions with topological information, creating a more refined and informative framework for settings where the underlying structure is crucial. This approach blends the uncertainty modeling of Gaussian priors, which provide a probabilistic framework with a known covariance structure, with topological insights that capture the shape, connectivity, and features of the parameter space. By integrating these components, we achieve a prior that reflects both the variability in the data and the inherent structural constraints of the system.

The proposed prior term is based on the concept of persistence distance $\|\cdot\|_{\mathrm{per}}$, and can be viewed as a generalization of the total variation (TV) prior. 
We see in Section 2 that for $1$-dimensional case, small values of $|q(v)-q(\tilde{v})|$ (or $q(x)-q(\tilde{x})$), which arise from closely paired extrema $(v, \tilde{v})$, are typically indicative of high-frequency oscillations or noise. On the other hand, large values of $|q(v)-q(\tilde{v})|$ correspond to prominent structural features of the function $q$. 
Let for simplicity 
\begin{align*}
P(q)=P_{+}\cup P_{-}\cup\{(v, \tilde{v})\}
\end{align*}
be the set of all (persistence) pairs, where $(v, \tilde{v})$ denotes the pair of knots whose corresponding function values are the global minimum and the global maximum of $q$. 
We define the prior as  
\begin{align}
\frac{{\rm d}\mu_{\rm pr}}{{\rm d}\mu_0}(q) \propto \exp\left(-\lambda\sum_{(x_j, \tilde{x}_j) \in P(q)} \alpha_j(q) |q(x_j) - q(\tilde{x}_j)|\right),
\end{align}
where $\lambda>0$ is a constant,  \( \alpha_j = \alpha_j({q}) = \alpha({q}, x_j, \tilde{x}_j) \) depends on the persistence \( |{q}(x_j) - q(\tilde{x}_j)| \).  
This corresponds to taking the regularization term \( R_{1\rm d}(q) \) as  
\begin{align}\label{eqn19}
R_{1\rm d}(q) = \lambda \sum_{(x_j, \tilde{x}_j) \in P(q)} \alpha_j(q) |q(x_j) - q(\tilde{x}_j)|.
\end{align}
The parameter \( \alpha_j \) should be large for small distances \( |q(x_j) - q(\tilde{x}_j)| \), emphasizing the penalization of minor oscillations or noise. Conversely, \( \alpha_j \) should be relatively small for large distances \( |q(x_j) - q(\tilde{x}_j)| \), allowing significant features of the function \(q\) to be preserved.
In \cite{Zheng2015}, the weight strategy is proposed to deal with denoising problems
\begin{align}\label{strategy}
\alpha_j(q)=(\kappa_j+1)\theta\frac{1}{1+\beta |q(\tilde{x}_j)-q(x_j)|}, 
\end{align}
where $\beta>0$ and $\kappa_j=\kappa(x_j, \tilde{x}_j)$ is the order of the pair $(x_j, \tilde{x}_j)$ in its chain of pairs and $\theta>1$. If we take $\alpha_j\equiv 1$ for all $j$, the proposed prior reduces to the TV prior.   In this sense, we say that it is a generalization of the total variation (TV) prior.

The concept of the $1$-dimensional persistence prior can be extended to the $2$-dimensional case. We only consider the case where $\Omega$ is a square domain. In this case, $q$ can be viewed as an image. The values of $q$ on vertices of the partition are arranged in a matrix.
We denote this matrix as $\mathrm Q=[\mathrm Q_{ij}]$, $i=1, \cdots, \mathrm I$, $j=1, \cdots, \mathrm J$, the $i$-th row as $\mathrm Q^{i}$, and the $j$-th column as $\mathrm Q_{j}$. The rows and columns of $\mathrm Q$ are then processed in a manner similar to the $1$-dimensional case. We treat each row (resp. column) one by one separately. 
To simplify the notation, we use $q^i$ and $q_j$ to denote the $i$-th row and the $j$-column respectively. 
 Then we impose a PH prior on each row (resp. column) and sum them, i.e., we take
 \begin{align}
 R_{2\rm d}(q)=\lambda\left(\sum_{i=1}^{\mathrm I}R_{1\rm d}({q}^{i})+\sum_{j=1}^{\mathrm J}R_{1\rm d}({q}_{j})\right), \,\, \lambda>0
 \end{align}
 and set
\begin{align}
\frac{{\rm d}\mu_{\rm pr}}{{\rm d}\mu_0}(q) \propto \exp\left(-R_{2\rm d}(q)\right).
\end{align}
This approach is equivalent to monitoring the topological structure of slices of the two-dimensional function along each direction. Although it cannot fully capture the topological information of the entire 2d function, it still provides valuable insights by preserving key topological features along specific slices.


If we denote both $R_{1d}$ and $R_{2d}$ by $R_{\rm PH}$, the posterior can be expressed as
\begin{align}\label{eqnph}
\frac{d\mu_{\rm post}}{d\mu_0} = \exp\left(-\Phi(q; u^\eta) - R_{\rm PH}(q)\right).
\end{align}
This formulation, together with \eqref{eqntv_c}, adopts the general form of \eqref{eqn18}:
\begin{align}\label{eqn24}
\frac{d\mu_{\rm post}}{d\mu_0}(q) = \exp(-V(q)) = \exp(-\Phi(q; u) - R(q)),
\end{align}
where $R$ represents either $R_{\rm PH}$ or $R_{\rm TV}$.

In this framework, $V(q)$ serves as the posterior potential functional. It consists of the data misfit term $\Phi(q; u^\eta)$ combined with the regularization term $R_{\rm PH}(q)$ (or $R_{\rm TV}(q)$). This potential functional quantifies the adjustment relative to the Gaussian reference measure $\mu_0$ that defines the posterior distribution.

\section{Delayed Acceptance Metropolis–Hastings Algorithm}

We rewrite \eqref{eqn24} as
\begin{align}\label{evn25}
\begin{aligned}
\frac{d\mu_{\rm post}}{d\mu_0}(q)=\exp(-V(q))&=\exp(-\Phi(q; u)-(1-a)R(q)-a R(q))\\
&:=\exp(-\varphi-a R(q)),
\end{aligned}
\end{align}
where $0\leq a\leq 1$. 
We design a two-stage acceptance Metropolis-Hastings (MH) algorithm based on the factorization. 
The decomposition naturally leads to a two-step acceptance mechanism:
\begin{enumerate}[i)]
\item First Stage (Pre-screening): Use the pCN proposal to generate a candidate $q'$. Accept $u'$ with probability 
\begin{align}\label{prob26}
\upsilon_1(q, q')=\min\left\{1, \exp(-a(R(q')-R(q)))\right\}.
\end{align} 
This stage filters proposals according to the partial weight $\exp(-aR)$.
\item Second Stage (Correction): If the proposal passes the first stage, apply a second acceptance/rejection step with probability
\begin{align}
\upsilon_2(q, q')=\min\left\{1, \exp\left(-\varphi(q')+\varphi(q)\right)\right\}.
\end{align}
This correction ensures that the overall algorithm targets the full posterior distribution. 
\end{enumerate}

By construction, the product of the two acceptance probabilities corresponds to the standard MH acceptance probability:
\begin{align}
\upsilon(q, q')=\upsilon_1(q, q')\cdot \upsilon_2(q, q')=\min\left\{1, \exp\left(-[\Phi(q')+R(q')]+[\Phi(q)+R(q)]\right)\right\},
\end{align}
which satisfies detailed balance with respect to $\mu_{\rm post}$. This scheme can be interpreted as a two-stage acceptance method, where the first step acts as a computationally cheaper pre-screening based on $R$,
 while the second step enforces correctness with respect to the full posterior. This approach is particularly useful when evaluation of $\Phi$ is computationally expensive, but $R$ is inexpensive to compute. Furthermore, for all values of $a$ in the interval $0<a\leq 1$, the computational cost 
 and effectiveness are identical. This implies that it is sufficient to consider $a=0$ or $a=1$. When $a=0$, the algorithm reduces to the standard pCN scheme.

\begin{algorithm}
\begin{algorithmic}[1]
\State \textbf{Initialization:}  Set sample number $N$ and $q^{(0)}=0$. Compute the corresponding negative log likelihood function $\Phi(q^{(0)}; u^\eta)$ and regularization term $R(q^{(0)})$. 
\State Move the sample to a proposal $\hat{q}^{(n)}=\sqrt{1-\rho^2}q^{(n-1)}+\rho \xi$, where $\xi\sim \mu_0$. 
\State Accept $\hat{q}^{(n)}$ with probability $\upsilon_1(q^{(n)}, \hat{q}^{(n)})$. 
\State If $\hat{q}^{(n)}$ is accepted in the last step, we accept this sample with probability $\upsilon_2(q^{(n)}, \hat{q}^{(n)})$. Otherwise, return step 2. 
\State When $n<N$, implement step 2 - step 4.
\end{algorithmic}
\caption{Delayed acceptance MH algorithm}
\label{alg3}
\end{algorithm}

\section{Numerical tests}

In this section, we present numerical experiments to validate the effectiveness of the proposed prior. We compare the performance of three different priors: the standard Gaussian prior, the TV-Gaussian prior, and our proposed PH-Gaussian prior. The results demonstrate the advantages of our approach in terms of accuracy and robustness.

 In all numerical tests,
the data is generated by adding a $1\text{\textperthousand}$ relative error to the numerical solution using the exact $q$, i.e.,
\begin{align}\label{dta}
u^\eta(x_o)=u(q)(x_o)+\eta(x_o), \quad \eta=1\text{\textperthousand}\|u(q)\|\xi,
\end{align} 
where $u(q)$ is the numerical solution and $\xi\sim N(0, I)$ and $x_o$ is observational points uniformly distributed in $\Omega$. And the error between the $n$-th samples and the exact solution $q$ is computed by
\begin{align}\label{errorformu}
\|q^{(n)}-q\|=\frac{1}{\sqrt{m-1}}\left(\sum_{i=1}^m |q^{(n)}(x_i)-q(x_i)|^2\right)^{\frac{1}{2}},
\end{align}
where $\{x_i\}$ are the space discretized mesh points. 
We employ the pCN algorithm (see Algorithm \ref{alg3}) to generate posterior samples from the target distribution \cite{Yao2016}. To ensure convergence to the stationary distribution and mitigate the effects of initial transients, we incorporate a burn-in phase and apply thinning (lag) to reduce sample autocorrelation. Specifically, we discard the first half of the generated samples as burn-in and set a lag of 
$5$ between retained samples. The sample mean is then computed from the remaining samples to estimate the parameter 
$q$.

\subsection{1d case}

We consider several different $1$d potential functions and generate $N=10^5$ samples to verify our method. 

\begin{example}\label{example5} First, we consider a function that is derived from the Weierstrass function
\begin{align*}
W(x)=\sum_{n=0}^\infty a^n \cos(b^n \pi x),
\end{align*}
where $0<a<1$, $b$ is a positive odd integer and $ab>1+\frac{3\pi}{2}$.  By $W_K(x)$ we denote the truncated series to the first $K+1$ terms. Take the exact $q$ as $q(x)=\frac{2}{\pi}\arctan(W_K(x))+1$ with $a=0.4$, $b=4$ and $K=10$. Here, we do not adhere to the strict parameter rule in the Weierstrass function.
\end{example}


In the sampling process for this example, we apply a transform $q(x)=\exp(g(x))$ and treat $g$ as the unknown target. The mean squared exponential kernel Gaussian is used as the base measure. We set $l=0.01$ and the weight parameter $\lambda=5$.
 The Gaussian prior produces  oscillatory results, whereas both the TV-Gaussian and PH-Gaussian priors achieve robust performance (see Fig. \ref{weierstrass}). We compare the numerical results for $a=0$ and $a=1$ for TV-Gaussian and PH-Gaussian priors using Algorithm \ref{alg3}. The  results demonstrate that the performance is almost identical for both $a=0$ and $a=1$ for each prior. However the time cost differs between two cases (see Table \ref{table_a}).  We observe that the delayed acceptance algorithm is slightly faster than the standard pCN scheme.  While the PH-Gaussian prior is computationally more expensive, it achieves higher accuracy in numerical reconstruction compared to the TV-Gaussian prior. For simplicity, we only consider the case $a=1$ in what follows. 
 \begin{table}[htbp]
 \begin{tabular}{c|c|c}
a& TV-Gaussian & PH-Gaussian \\ 
\hline
$0$&$14.51$&$3144.74$ \\
\hline
$1$&$11.75$& $3123.19$ 
\end{tabular}
\caption{Sampling Time (s) for Example \ref{example5}}
\label{table_a}
\end{table}
\begin{figure}[h] %
    \centering
    \includegraphics[width=0.32\textwidth]{weirstrass_gauss.png} %
    \includegraphics[width=0.32\textwidth]{weirstrass_tv.png} %
     \includegraphics[width=0.32\textwidth]{weirstrass_ph.png} %
    \caption{Numerical comparisons using different priors for Example  \ref{example5}: Left: Gaussian prior; Middle: TV-Gaussian  prior; Right: PH-Gaussian prior.}
    \label{weierstrass}
\end{figure}

\begin{figure}[htbp]
    \centering
    
    \subfigure[Example \ref{example1}]{
        \includegraphics[width=0.4\textwidth]{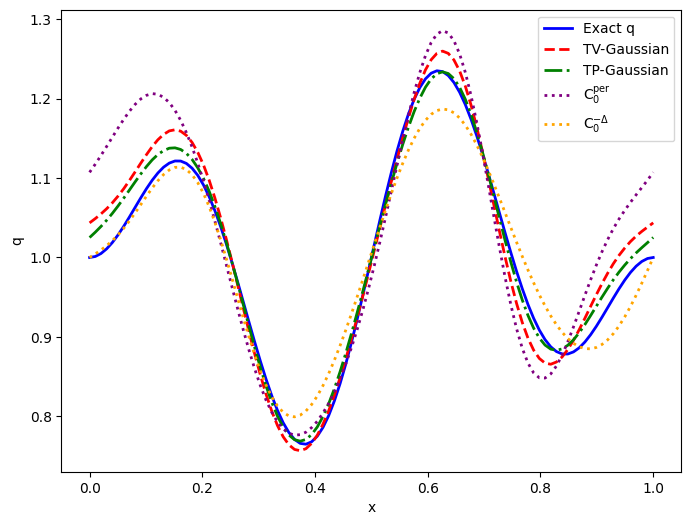}
        \label{subfig1}
    }
    \hfill
    \subfigure[Example \ref{example2}]{
        \includegraphics[width=0.4\textwidth]{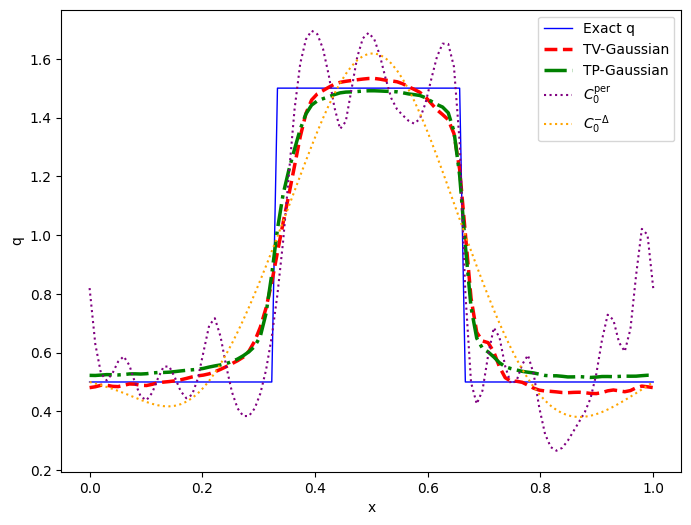}
        \label{subfig2}
    }

    \subfigure[Example \ref{example3}]{
        \includegraphics[width=0.4\textwidth]{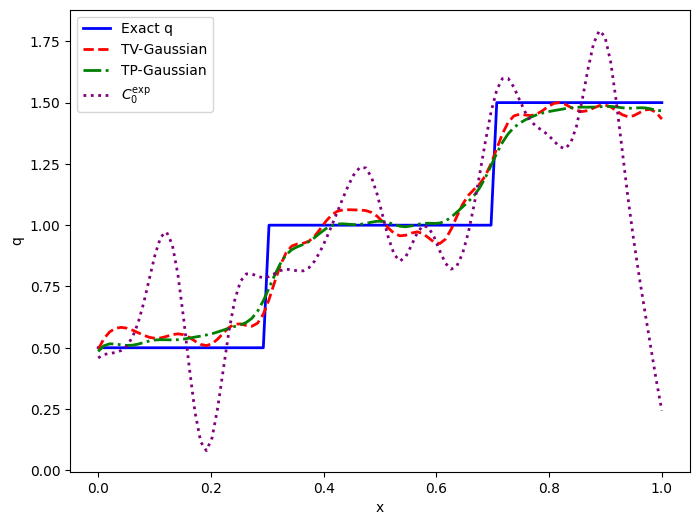}
        \label{subfig3}
    }
    \hfill
    \subfigure[Example \ref{example4}]{
        \includegraphics[width=0.4\textwidth]{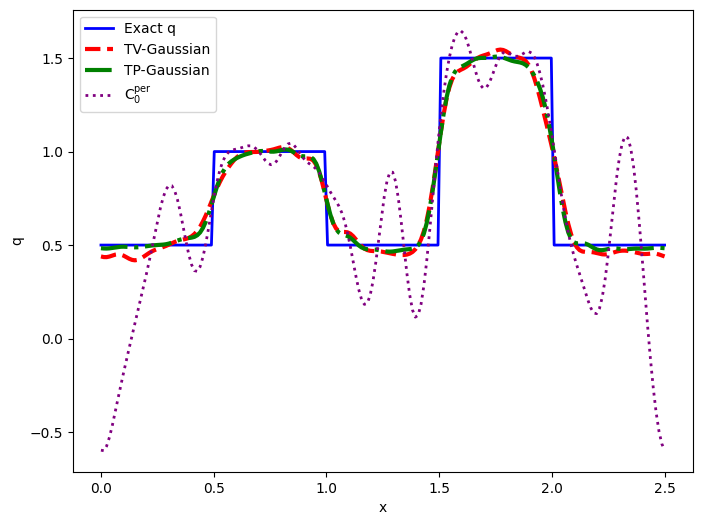}
        \label{subfig4}
    }
    
    \caption{Numerical comparisons using different priors for 1d case.}
    \label{test1}
\end{figure}

\begin{example}\label{example1} In this example, we consider a smooth 1-dimensional function in $\Omega=(0, 1)$
\begin{align*}
q(x)=1+x(1-x)\sin(4\pi x).
\end{align*}
\end{example}
In this experiment, we compare the effectiveness of different priors introduced in Section 3: the Gaussian priors with $C_0^{-\Delta}$ (with $s=1.5$) and the periodic mean exponential kernel $C_0^{\rm per}$, as well as the TV-Gaussian prior and the PH-Gaussian prior.
For the latter two priors, we use the periodic mean squared exponential kernel Gaussian in $\mu_0$ as the reference measure. The parameters for these priors are set as follows: $l = 1.0$, $p = 1$ for $C_0^{\rm per}$, and $\theta = 3$, $\beta = 0.001$ in the PH prior as given in \eqref{strategy}.
The weight parameter $\lambda$ in both the TV-Gaussian and PH-Gaussian priors is set to 4. The proposal parameter $\rho$ in Algorithm \ref{alg3} is set to 0.002. 
Fig.~\ref{test1}--\subref{subfig1} presents numerical reconstructions that illustrate the efficacy of the proposed prior for smooth scenarios. Notably, the PH-Gaussian prior demonstrates competitive performance compared to alternative priors.


\begin{example}\label{example2}
In the third example, we consider the exact \( q \) as a 1-dimensional piecewise constant function on $\Omega=(0, 1)$:
\[
q(x) = \begin{cases}
1.5, & \frac{1}{3} \leq x < \frac{2}{3}, \\
0.5, & \text{otherwise}.
\end{cases}
\]
\end{example}
This numerical experiment demonstrates the performance of both the baseline prior from Example \ref{example1} and the proposed prior $\mu_0$ with the kernel $C_0^{-\Delta}$ (with parameter $s=1.6$). 
For consistency with Example \ref{example1}, we explicitly list only the modified parameters: the length scale $l=0.2$ in $C_0^{\rm per}$, the regularization weight $\lambda=20$, and the proposal step size $\rho=0.01$.
The numerical comparison is given in Fig.~\ref{test1}--\subref{subfig2}. Obviously, when dealing with a function that has jumps, using only a Gaussian prior leads to unsatisfactory results. Specifically, the Gaussian prior \( \mu_0 \) with \( C_0^{-\Delta} \) enforces an excessively high degree of smoothness on the function due to its highly smooth eigenfunctions within the domain. This makes it difficult to capture local features, especially for non-smooth functions. Whether using the Karhunen-Loève expansion, finite difference methods, or finite element methods to discretize the covariance operator \( C_0^{-\Delta} \), we encounter the problem that either the number of expansion terms becomes too large, or the discretization grid is limited by machine precision, preventing the capture of local features. 
When using \( \mu_0 \) with \( C_0^{\rm per} \), the result exhibits oscillations, likely arising from the properties of the periodic exponential covariance kernel, particularly its inherent smoothness and long-range dependencies, which can introduce oscillatory artifacts into the solution.
In contrast, the TV-Gaussian and PH-Gaussian priors can significantly alleviate these issues. In the subsequent examples, we will no longer employ \( \mu_0 \) with the covariance operator \( C_0^{-\Delta} \), as our focus is solely on testing non-smooth scenarios.


\begin{figure}[htbp]
    \centering

    \begin{minipage}[b]{0.3\textwidth}
        \includegraphics[width=\textwidth]{weierstrass_trace.png}
        \caption*{Example  \ref{example5}}
        \label{fig:subfig5}
    \end{minipage}
    \begin{minipage}[b]{0.3\textwidth}
        \includegraphics[width=\textwidth]{ex1_trace_plot.png}
        \caption*{Example  \ref{example1}}
        \label{fig:subfig1}
    \end{minipage}
    \begin{minipage}[b]{0.3\textwidth}
        \includegraphics[width=\textwidth]{ex2_trace_plot.png}
        \caption*{Example \ref{example2}}
        \label{fig:subfig2}
    \end{minipage}
    \begin{minipage}[b]{0.3\textwidth}
        \includegraphics[width=\textwidth]{ex3_trace_plot.png}
        \caption*{Example \ref{example3}}
        \label{fig:subfig3}
    \end{minipage}
    \begin{minipage}[b]{0.3\textwidth}
        \includegraphics[width=\textwidth]{ex4_trace_plot.png}
        \caption*{Example \ref{example4}}
         \label{fig:subfig4}
    \end{minipage}
   
    \caption{The error traces for Example \ref{example5}-\ref{example4}.}
    \label{trace_plot}
\end{figure}

\begin{example}\label{example3} We test the proposed prior on a 1D step function, 
\begin{align*}
q(x)=\begin{cases}
0.5, & 0 \leq x < 0.3, \\
1, &0.3\leq x<0.7,\\
1.5, & 0.7\leq x\leq 1.
\end{cases}
\end{align*}
\end{example}
In this example, the parameter \( l \) in the Gaussian measure $\mu_0$ is set to \( 0.05 \) in all cases. The weight parameter \( \lambda \) is set to \( 10 \), and the Markov transition parameter \( \beta \) is set to \( 0.005 \). Other parameters are chosen to be the same as in Example  \ref{example2}. The numerical results are displayed in Fig.~\ref{test1}--\subref{subfig3}. As in Example  \ref{example2}, the oscillation occurs when using the prior $\mu_0$ with the exponential covariance kernel and the reason is similar to that in Example  \ref{example2} . In addition, the TV-Gaussian and PH-Gaussian priors can fulfill the case well. 

\begin{example}\label{example4} We adjust the piecewise constant in Example  \ref{example2}  to a more complicated case as 
\begin{align*}
q(x)=\begin{cases}
0.5, & 0 \leq x < 0.5, \\
1, &0.5\leq x<1.0,\\
0.5, & 1.0\leq x\leq 1.5,\\
1.5, & 1.5\leq x<2.0,\\
0.5, & 2.0\leq x\leq 2.5.
\end{cases}
\end{align*}
\end{example}
We  use the periodic mean squared exponential kernel Gauss as the base measure with $l=0.2$, $p=2.5$.   The weight parameter $\lambda=20$. The parameter of the PH-Gaussian prior  $\theta=3$ and $\beta=0.001$ in \eqref{strategy}. 
From the displayed results (see Fig.~\ref{test1}--\subref{subfig4}), it is evident that the same phenomenon observed in Examples  \ref{example2}  and  \ref{example3}  is replicated here.

We plot the error traces  for Examples  \ref{example5}-\ref{example4}   in Fig. \ref{trace_plot}, with the error computed according to \eqref{errorformu}. These results demonstrate that the proposed prior delivers strong numerical performance, particularly for functions with discontinuities. While the PH and TV priors yield comparable results in most cases, a discernible difference emerges in Example  \ref{example5}  when applied to a highly oscillatory function—here, the advantage of our method becomes visually apparent, as reflected in the error trace plot.

\subsection{2d case}
 \begin{example}\label{example6}  We consider a 2d piecewise constant function in a squared domain $\Omega=(0, 1)\times (0, 1)$
\begin{align*}
q(x, y)=\left\{
\begin{aligned}
&1.5,&& (x-0.5)^2+(y-0.5)^2\leq 0.25^2, \\
&0.5, & &\text{otherwise}.
\end{aligned}
\right.
\end{align*}
\end{example}
In this example, we take \( \mu_0 = N(0, C_0^{\exp}) \) with the parameter \( l = 0.1 \). We generate \( N = 40,000 \) samples.  
The weight parameter \( \lambda \) in the TV-Gaussian and PH-Gaussian priors is set to 2. Other parameters remain the same as those in Example  \ref{example1}. 
\begin{figure}[ht]
\centering
\subfigure[$C_0^{\exp}$]{
\includegraphics[width=0.3\textwidth]{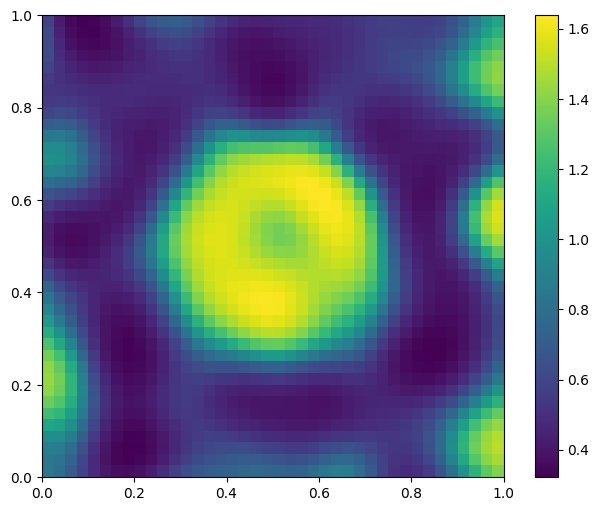}}
\subfigure[TV-Gaussian prior]{
\includegraphics[width=0.3\textwidth]{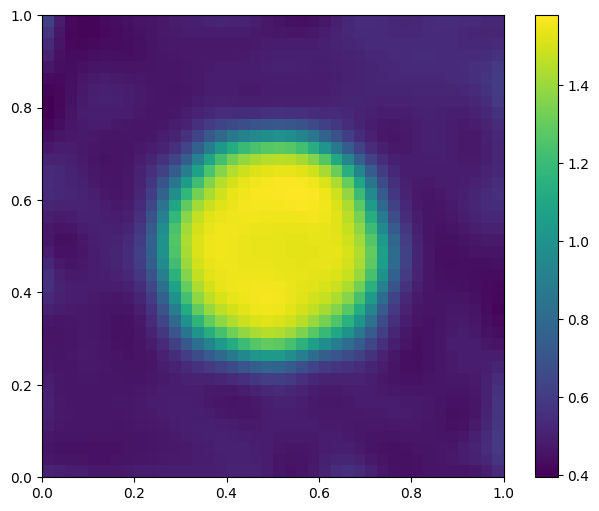}
}
\subfigure[PH-Gaussian prior]{
\includegraphics[width=0.3\textwidth]{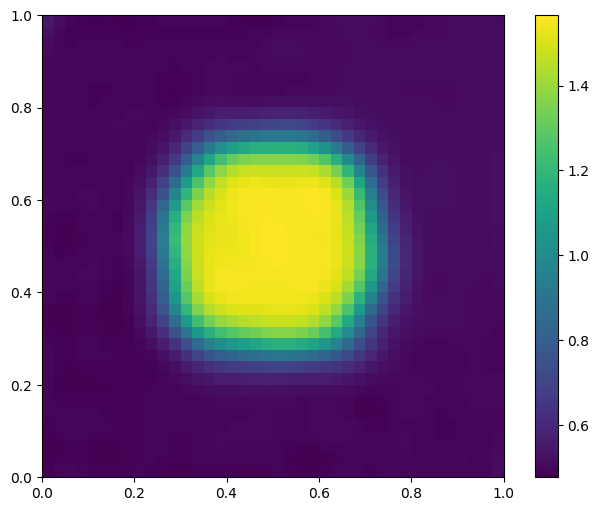}
}
\subfigure[The error trace plot]{
\includegraphics[width=0.4\textwidth]{ex5_trace_plot.png}
}
\caption{Numerical illustration using different priors for Example  \ref{example6}: (a) the Gaussian prior with $C_0^{\exp}$; (b) the TV-Gaussian prior; (c) the PH-Gaussian prior; (d) the error trace. }
\label{2D_1}
\end{figure}

\begin{example}\label{example7}We consider a 2d piecewise constant function in a squared domain $\Omega=(0, 2)\times (0, 2)$
\begin{align*}
q(x, y)=\left\{
\begin{aligned}
&1,&& (x-0.6)^2+(y-0.6)^2\leq 0.3^2, \\
&1.5, && (x-1.4)^2+(y-1.4)^2\leq 0.3^2\\
&0.5, & &\text{otherwise}.
\end{aligned}
\right.
\end{align*}
\end{example}

\begin{figure}[ht]
\centering
\subfigure[$C_0^{\exp}$]{
\includegraphics[width=0.3\textwidth]{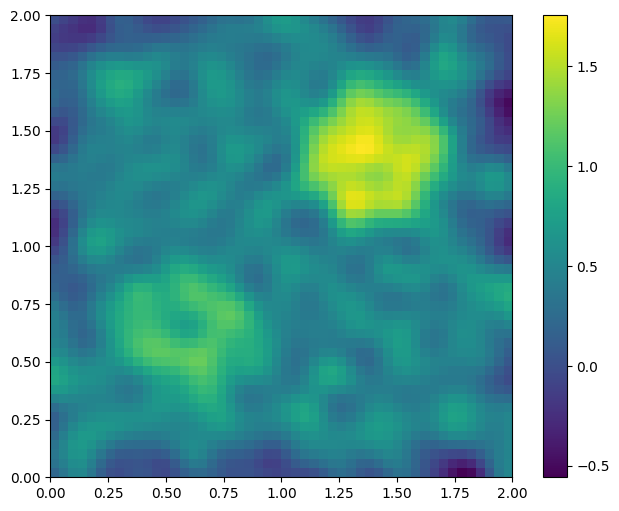}}
\subfigure[TV-Gaussian prior]{
\includegraphics[width=0.3\textwidth]{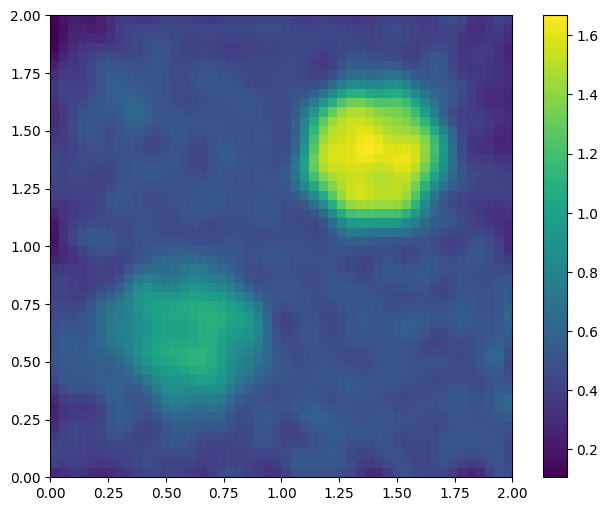}
}
\subfigure[PH-Gaussian prior]{
\includegraphics[width=0.3\textwidth]{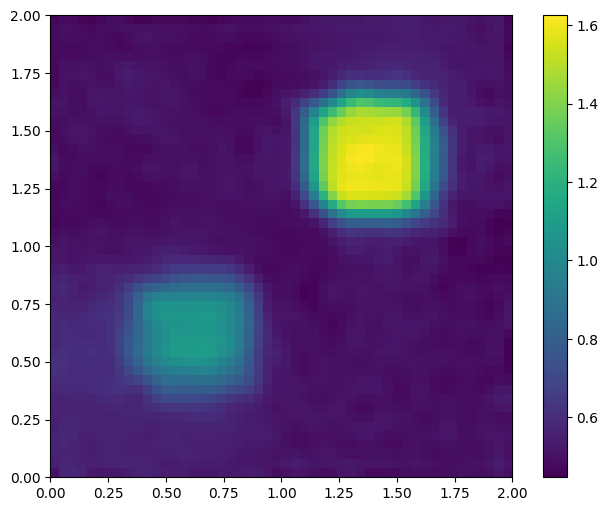}
}
\subfigure[The error trace plot]{
\includegraphics[width=0.4\textwidth]{ex6_trace_plot.png}
}
\caption{Numerical illustration using different priors for Example  \ref{example7}: (a) the Gaussian prior with $C_0^{\exp}$; (b) the TV-Gaussian prior; (c) the PH-Gaussian prior; (d) the error trace. }
\label{2D_2}
\end{figure}
In this example, apart from adjusting the parameter $l$ to $0.1$ and $\lambda$ to $5$, all other parameters remain the same as in Example  \ref{example5}. We generate $N=80,000$ samples.

From the results displayed in Fig.~\ref{2D_1}-\ref{2D_2} for  2D examples, we observe conclusions nearly identical to those in the 1D case. This confirms that the PH-Gaussian prior yields reliable reconstruction results, performing competitively with the TV-Gaussian prior and outperforming the standard Gaussian prior alone. The error trace plots further demonstrate that the PH-Gaussian prior exhibits superior numerical effectiveness compared to the TV-Gaussian prior when reconstructing functions with complex topological structures, as evidenced in Example  \ref{example7}.

\section{Conclusion}
We employ persistent homology as a tool to construct a hybrid prior for estimating the unknown variable in the inverse potential problem within the Bayesian framework. A key feature of this prior is that it constrains the topological variation of the unknown variable, effectively imposing a regularity condition akin to the total variation (TV). Moreover, the PH-based prior is defined on a topological space, which is a mild limitation. This flexibility allows the method to be easily extended to a broader range of applications.
In numerical practice, the PH-Gaussian prior demonstrates excellent performance. Its ability to produce robust results is enhanced by the flexibility of its parameter settings, which allow for more adaptable and precise tuning.

\bibliography{sn-bibliography}

\end{document}